\documentclass[10pt]{article}
\usepackage{amsmath, amsfonts, amssymb, latexsym}
\usepackage[english]{babel}
\usepackage{hyperref}
\usepackage{color}
\usepackage{pdfpages}

\newtheorem{theorem}{Theorem}[section]
\newtheorem{prop}[theorem]{Proposition}
\newtheorem{cor}[theorem]{Corollary}

\newtheorem{lemma}[theorem]{Lemma}
\newtheorem{definition}[theorem]{Definition}

\newtheorem{remark}[theorem]{Remark}
\numberwithin{equation}{section}
\input{tcilatex}
\begin{document}

\title{Limit theorem and LIL for some additive functionals associated to fBm
and Riemann-Liouville process}
\date{\today }
\author{Mohamed Ait Ouahra, Abderrahim Aslimani, \and Mhamed
Eddahbi, Mohamed Mellouk }
\maketitle

\begin{abstract}
In this paper, we first establish a strong approximation version for the
first order limit theorem of some additive functionals related to two
non-Markovian Gaussian processes: the fractional Brownian motion (fBm) and the
Riemann-Liouville process. As application, we give the law of iterated
logarithm of the same additive functionals associated to this two processes.%
\newline

\noindent {\footnotesize {\emph{Key words:} Limit theorem; Fractional
Brownian motion (fBm); Riemann-Liouville process; Additive functional; Law
of Iterated Logarithm (LIL); Local time; Strong approximation.} 
}
\end{abstract}



\section{Introduction.}

The first order limit theorem was first established by Darling and Kac \cite%
{DK} and affirms that the family of processes,
\begin{equation*}
\left\{ \frac{1}{\lambda ^{\frac{1}{2}}}\int_{0}^{\lambda t}f(W(s))ds,\
t\geq 0\right\}, \ \lambda >0,
\end{equation*}%
converges in law, in the space of continuous functions, as $\lambda $ goes
to infinity, to the process $\left\{ \overline{f}\ell (0,t),\ t\geq
0\right\} $, where $W=\{W(t),\ t\geq 0\}$ is a Brownian motion, $\ell (x,t)$
its local time at time $t\geq 0$ and level $x\in \mathbb{R}$ and
\begin{equation*}
\overline{f}=\int_{\mathbb{R}}f(x)dx\neq 0.
\end{equation*}%
A strong approximation result of the family of the processes
\begin{equation*}
\left\{ \int_{0}^{t}f(W(s))ds,\ t\geq 0\right\} ,
\end{equation*}%
was given by Csaki \emph{et al}. \cite{CCFR}: for all sufficiently small $%
\varepsilon >0$, when $t$ goes to infinity, we have
\begin{equation*}
\int_{0}^{t}f(W(s))ds=\overline{f}\ell (0,t)+o\left( t^{\frac{1}{2}%
-\varepsilon }\right) ,\ \text{\ a.s. }
\end{equation*}%
with $\int_{\mathbb{R}}|x|^{k}f(x)dx<+\infty $ for some $k>0$.


A strong approximation version for the first-order limit theorem of the
process
\begin{equation*}
\left\{ \int_{0}^{t}f(X(s))ds,\ t\geq 0\right\} ,
\end{equation*}%
associated to the symmetric stable process $\{X(t),\ t\geq 0\}$ of index $%
1<\alpha \leq 2$ was given by Ait Ouahra and Sghir \cite{AS}. This
generalizes the result of Csaki \emph{et al}. \cite{CCFR} for Brownian
motion which leads to give the law of iterated logarithm for some additive
functionals of this process.\newline
Kasahara and Kosugi, established in \cite{KK} a limit theorem version of
fractional Brownian motion $B^{H}=\{B^{H}(t),t\geq 0\}$ given as follows: if
$f\in L^{1}(\mathbb{R})$ such that $\overline{f}\neq 0$, then the family of
processes,
\begin{equation*}
\left\{ \frac{1}{\lambda ^{1-H}}\int_{0}^{\lambda t}f(B^{H}(s))ds,\ t\geq
0\right\}, \ 0<H<1,
\end{equation*}%
converges in law, in the space of continuous functions, as $\lambda $ goes
to infinity, to the process $\left\{ \overline{f}L^{H}(0,t),\ t\geq
0\right\} $, where $L^{H}(0,t)$ is the local time of $B^{H}$ at $x=0$.%
\newline
In the case of Riemann-Liouville process $W^{\beta }$ of parameter $\beta >0$%
, we can prove easily by using the dominated convergence theorem and the
fact that the local time has compact support, that for $f\in L^{1}(\mathbb{R}%
)$ such that $\overline{f}\neq 0$,
\begin{equation*}
\left\{ \frac{1}{\lambda ^{1-\beta }}\int_{0}^{\lambda t}f(W^{\beta
}(s))ds,\ t\geq 0\right\} ,
\end{equation*}%
converges in law, in the space of continuous functions, as $\lambda $ goes
to infinity, to the process $\left\{ \overline{f}L^{\beta }(0,t),\ t\geq
0\right\} $, where $L^{\beta }(0,t)$ denote the Riemann-Liouville local time
at $x=0$ when $0<\beta <1$. Throughout this paper, we use the same symbol $%
X^{\tau }=\left\{ X^{\tau }(t),\ t\geq 0\right\} $ to denote each of
Gaussian $\tau $-self-similar processes: The fractional Brownian motion ($%
\tau =H$, fBm for short) and the Riemann-Liouville process ($\tau =\beta >0$%
), and we denote $\left\{ L^{\tau }(x,t),\ t\geq 0,x\in \mathbb{R}\right\} $
its local time (under the condition of the existence $\tau <1$ for the
Riemann-Liouville process). \newline
We can then ask ourselves if the process $X^{\tau }$ possesses the strong
approximation version with respect to the same family of the processes
\begin{equation*}
\left\{ \int_{0}^{t}f(X^{\tau }(s))ds,\ t\geq 0\right\} ,
\end{equation*}%
in the sense that: for all sufficiently small $\varepsilon >0$, when $t$
goes to infinity, we have (under some conditions on $f$)
\begin{equation*}
\int_{0}^{t}f(X^{\tau }(s))ds=\overline{f}L^{\tau }(0,t)+o\left( t^{1-\tau
-\varepsilon }\right) ,\ \text{\ a.s. }
\end{equation*}%
In the classical case of Brownian motion and in the general case of the
symmetric stable processes, this property is a consequence of the spatial
homogeneity which follows from the stationarity of increments. However, the
process $X^{\tau }$ is not a Markov process but we show that we have an
almost similar property of some stochastic processes associated to the local
time $L^{\tau }$ (see lemma \ref{lem2.21} and the proof of lemma \ref{lem2.3}
).

\noindent \textbf{Notations and definitions.} We give here some notations
which can serve us throughout this paper: The underlying parameter space is \mbox{$%
\mathbb{R}^{+}=[0,+\infty)$} and we denote $\mathcal{B}(\mathbb{R}^{+})$, the
family of Borel subsets of $\mathbb{R}^{+}$. We shall write $dx$ (resp. $ds$%
) for Lebesgue's measure on $\mathbb{R}$ (resp. on $\mathbb{R}^+$). For
every Borel function $f$ on $\mathbb{R}$, we denote by $\overline{f}$ the
integral
\begin{equation*}
\overline{f}:=\int_{\mathbb{R}}f(x)dx.
\end{equation*}
Let $(\Omega,\mathcal{F},\mathbb{P})$ be a probability space, $X=\{X(t),t\geq0\}$ and $Y=\{Y(t), t\geq 0\}$ are two processes on this space.
We will use $X\overset{(d)}{=}Y$ to mean that the two processes $X$ and $Y$
have the same finite-dimensional distributions. The notation ``a.s." is the
abbreviated form of ``almost surely". If $X=\{X(t), \ t\geq0\}$ is a
stochastic process starting from $x\in\mathbb{R}$, a probability measure $%
\mathbb{P}^x $ is defined to be the distribution of the random variable $t\mapsto
X_t(\cdot)$, i.e.,
\begin{equation*}
\mathbb{P}^x(A):=\mathbb{P}\left((t\mapsto X_t)\in A\right), \ \ A\in\mathcal{F}.
\end{equation*}
The expectation under $\mathbb{P}^x$ is denoted by $\mathbb{E}_x$. In addition, the
norm $\left(\mathbb{E}_0|\cdot|^p\right)^{1/p}$ of the space $L^{p}(\Omega)$ with exponent $1\leq p<+\infty$, will be denoted by $%
\left
\|\cdot\right\|_p$.\newline
Finally, we will also use in this paper, unspecified positive finite
constants $C_1$, $C_2,\dots$, which may not necessarily be the same in each
occurrence.

The rest of this paper is organized as follows: in the next section, we
present some basic facts about the processes studied in this paper and a
brief reminder on the local time of Gaussian processes. In section 3, we
present our main result about the strong approximation of the additive
functionals associated to the process $X^{\tau}$
\begin{equation*}
\left\{\int_{0}^tf(X^{\tau}(s))ds , \ t\geq0\right\},
\end{equation*}
with $\overline{f}:=\int_{\mathbb{R}}f(x)dx\neq0$ and $\int_{\mathbb{R}%
}|x|^kf(x)dx<+\infty$ for some $k>0$. As application, we give in the last
section the law of iterated logarithm of this additive functionals
associated to this two Gaussian processes.

\section{Essential concepts and results about fBm, Riemann-Liouville process
and Local time of Gaussian processes}


\subsection{On the fractional Brownian motion (fBm)}

The fractional Brownian motion $B^{H}$ with Hurst exponent $%
H\in (0,1)$ introduced by Kolmogorov \cite{Kol} in a study of turbulence
have attracted researchers in recent years. By definition, it is an a.s.
continuous zero-centered Gaussian process with its covariance function given
by
\begin{equation*}
R(t,s)=\mathbb{E}\left[ B^{H}(t)B^{H}(s)\right] =\frac{1}{2}\left(
t^{2H}+s^{2H}-|t-s|^{2H}\right) .
\end{equation*}%
His study was resumed and extended by Mandelbrot and van Ness \cite{MV}.
This process is a self-similar process, that is,
\begin{equation*}
\{B^{H}(at),\ t\geq 0\}\overset{(d)}{=}\{a^{H}B^{H}(t),\ t\geq 0\},
\end{equation*}%
for all $a>0$. This property shows that a change of scale in time is
equivalent (in law) to a change of scale in space.

The fBm is also a process with stationary increments, that is,
\begin{equation*}
\left\{B^H(t+s)-B^H(s), \ t\geq0\right\} \overset{(d)}{= }%
\left\{B^H(t)-B^H(0), \ t\geq0\right\},
\end{equation*}
for all $s\in\mathbb{R}^+$. The self-similarity and stationarity of the
increments are two main properties for which fBm enjoyed success as a
modeling tool. In fact, fBm is the only Gaussian self-similar process with
stationary increments.

%
%


\begin{remark}
\label{rmk1.02} The fBm is neither a semi-martingale, nor a Markov process.
\end{remark}

\subsection{On the Riemann-Liouville process}

The Riemann-Liouville process $W^{\beta }(t)$ with index $\beta >0$ is
defined by the following stochastic integral, called moving average
representation which was introduced by Mandelbrot and Van Ness \cite{MV}
\begin{equation*}
W^{\beta }(t)=\int_{0}^{t}(t-s)^{\beta -1/2}dW(s),\ \ t\in \mathbb{R}^{+}
\end{equation*}%
where $W=\{W(t), \ t\geq 0\}$ is a standard Brownian motion.\newline
The process $W^{\beta }$ is a self similar zero-mean Gaussian process with
index $\beta $, but $W^{\beta }$ does not have stationary increments and
there is no upper bound restriction on index $\beta >0$.\newline
The relation between Riemann-Liouville process and fBm becomes transparent
when we write a moving average representation of $\{B^{H}(t), \ t\in \mathbb{%
R} \}$, in the form
\begin{equation*}
B^{H}(t)=c_{H}\int_{-\infty }^{t}\left[ (t-s)^{H-1/2}-(-s)_{+}^{H-1/2}\right]
dW(s),
\end{equation*}%
where $$
c_{H}=\sqrt{2H}2^{H}\beta (1-H,H+1/2)^{-1/2}
$$ and $$\beta (a,b)=\int_{0}^{1}x^{a-1}(1-x)^{b-1}dx$$ is the usual beta
function. Then we can deduce that
\begin{equation*}
B^{H}(t)=c_{H}(W^{H}(t)+Z^{H}(t)),
\end{equation*}%
where $Z^{H}$ is a process independent of $W^{H}$ given by
\begin{equation*}
Z^{H}(t)=\int_{-\infty }^{0}\left[ (t-s)^{H-1/2}-(-s)^{H-1/2}\right] dW(s).
\end{equation*}

\subsection{Local time of LND Gaussian process}

In this section we give a brief review of the concepts of local times
theory. For an excellent summary of local times results, the reader is
referred to Berman \cite{B}, Geman and Horowitz \cite{GH}, Xiao \cite{X},
see also Adler \cite{A}.

\begin{definition}
\label{def 1.1}Let $X = \{X(t), \ t\geq0\}$ be a real-valued separable
random process with Borel sample functions defined on a general probability
space $(\Omega,\mathcal{F},\mathbb{P})$. For any Borel set $B\subset \mathbb{%
R}^+$, the occupation measure of $X$ on $B$ is defined as
\begin{eqnarray}
\mu_B(A): = \lambda\left\{s\in B: \ X(s) \in A\right\}=\int_{B}\mathbb{I}%
_{A}(X(t))dt, \label{eq2}
\end{eqnarray}
for all $A \in \mathcal{B}(\mathbb{R})$, and where $\lambda$ is the
one-dimensional Lebesgue measure on $\mathbb{R}^+$. If $\mu_B$ is absolutely
continuous with respect to Lebesgue measure on $\mathbb{R}$, we say that $X$
has a local time on $B$ and define its local time $L( \cdot,B)$, to be the
Radon-Nikodym derivative of $\mu_B$. Here, $x$ is the so-called space
variable and $B$ is the time variable of the local times.
\end{definition}

Notice that if $X$ has local times on $\mathbb{R}^{+}$ then for every Borel
set $B\subset \mathbb{R}^{+}$, $L(x,B)$ also exists. Moreover, by standard
monotone class arguments, one can deduce that the local times have a
measurable modification that satisfies the following occupation density
formula: for every Borel set $B\in \mathcal{B}(\mathbb{R}^{+})$ and every
measurable function $f:\mathbb{R}\longrightarrow \mathbb{R}^{+}$,
\begin{eqnarray}
\int_B f(X(t))dt= \int_{\mathbb{R}} f(x)L(x,B)dx, \label{eq3}
\end{eqnarray}%
which holds for any Borel function $f$ and follows from the definition of
Radon-Nikodym derivative and linearity combined with a monotone convergence
argument. Moreover, from the occupation density formula, the local time
admits the following almost sure approximation
\begin{eqnarray}
L(x,t)=\lim_{\varepsilon \rightarrow 0^{+}}\frac{1}{2\varepsilon }%
\int_{0}^{t}\mathbb{I}_{\{\mid X(s)-x\mid \leq \varepsilon \}}ds.
\label{eq08}
\end{eqnarray}

\begin{remark}
\label{rmk1.1}For a fixed rectangle $D=\prod_{i=1}^{m}[a_{i},a_{i}+h_{i}]$, $0\leq a_i$, $0<h_{i}$ for $i=1,\ldots ,m$, if we can choose $L(x,\prod_{i=1}^{m}[a_{i},a_{i}+t_{i}])$, to be a
continuous function of $(x,t_{1},\ldots ,t_{m})$, $x\in \mathbb{R}$, $0\leq
t_{i}\leq h_{i}$ for $i=1,\ldots ,m$, then $X$ is said to have a jointly
continuous local time on $D$. Under this condition, $L(x,\cdot )$ can be
extended to be a finite measure (called Hausdorff measure) supported on the
level set $\{t\in D:X(t)=x\}$ (see \cite{A} for more details). In other
words, local times often act as a natural measure on the level sets of $X$.
\end{remark}

From now on, we write $L(x,t)$ instead of $L(x,[0, t ])$.

\begin{remark}
\label{rmk1.2}\ A necessary and sufficient condition of the existence for a
space-square integrable local time $L(x,t)$ of a real measurable process $%
X=\{X(t), \ t\in[0,T]\}$ for $T>0$ defined on a probability space $(\Omega,%
\mathcal{F},\mathbb{P})$ was established by Geman and Horowitz \cite{GH} and
given by
\begin{eqnarray}
\liminf_{\varepsilon\to 0^+}\frac{1}{\varepsilon} \int_{0}^{T}\int_{0}^{T}%
\mathbb{P}\left( \left|X(t)-X(s) \right|\leq \varepsilon \right)dsdt<+\infty.
\label{eq4}
\end{eqnarray}
\end{remark}

\begin{remark}
\label{rmk1.3}By standard monotone class argument, one can deduce that the
local times have a version, still denoted by $L$, such that it is a kernel
in the following sense:\newline
(i) For each fixed $B \in \mathcal{B}(\mathbb{R}^+)$, the function $x
\mapsto L(x, B)$ is Borel measurable in $x \in \mathbb{R}$. \newline
(ii) For every $x \in\mathbb{R} $, $B\mapsto L(x,B)$ is a Borel measure on $%
\mathcal{B}(\mathbb{R}^+)$.
\end{remark}

Berman \cite{B} developed Fourier analytic methods for studying the
existence and regularity of the local times of Gaussian processes. His
methods were extended by Pitt \cite{P} and Geman and Horowitz \cite{GH} to
Gaussian random fields.\newline
Let $X=\{X(t),\ t\geq 0\}$ be a Gaussian random field with values in $%
\mathbb{R}$. Let us give here an outline of the analytic method used by
Berman for the calculation of the moments of local times. For a fixed sample
function at fixed $t$, the Fourier transform on $x$ of $L(x,t)$ is the
function
\begin{equation*}
\widehat{L(x,t)}(u)=\int_{\mathbb{R}}\exp (iux)L(x,t)dx.
\end{equation*}%
Using the density of occupation formula we get
\begin{equation*}
\widehat{L(x,t)}(u)=\int_{0}^{t}\exp \left( iuX(s)\right) ds.
\end{equation*}%
Therefore, the local time may be represented as the inverse Fourier
transform of this function, i.e.
\begin{equation}
L(x,t)=\frac{1}{2\pi }\int_{\mathbb{R}}\left( \int_{0}^{t}\exp \left[
iu\left( X(s)-x\right) \right] ds\right) du. \label{eq7}
\end{equation}%
On the other hand, Berman \cite{B} and Geman and Horowitz \cite{GH} have
shown that, any continuous locally non deterministic Gaussian (LND in
abbreviated form) process $X$ has a jointly continuous version of local time
$(x,t)\mapsto L(x,t)$ with respect to time-space variables,
and if $X$ has a stationary increments and incremental variance $\sigma ^{2}(t)\backsim
t^{2\gamma }$ for small $t>0$ and $0<\gamma \leq 1$ then the mapping $t\mapsto L(x,t)$ is H\"{o}lder continuous of order $\gamma ^{\prime }$ for every $\gamma ^{\prime }<1-\gamma $ uniformly in $x$ from any compact set, and the map $x\mapsto L(x,t)$ is H\"{o}lder continuous of order $\gamma
^{\prime \prime }$ for every $\gamma^{\prime \prime }<\frac{1}{2}(\frac{1}{\gamma }-1)$ uniformly in $t$ on every compact interval.\\ In particular, for
any $0<\tau <1$, there exists a version of local time $L^{\tau }(x,t)$ of
the process $X^{\tau }$ such that the map $L^{\tau }(x,t)$ is a.s.
continuous in both $t$ and $x$. Moreover, it can be chosen so that the map $%
x\mapsto L^{\tau }(x,t)$ is H\"{o}lder continuous of order $\gamma_{1}<\frac{1}{2}(\frac{1}{\tau }-1)$ as well as uniformly continuous in $t$ on every compact interval. Moreover, according to Xiao \cite{X}, we have the following result:

\begin{lemma}
\label{LemReg} Let $T>0$. There exist constants $C_1>0$ and $C_2>0$, such that for any $p\geq 1$, $0\leq t,s\leq T$ and $x,y\in \mathbb{R}$ we have
\begin{eqnarray} \label{reg1}
\|L^{\tau}(x,t)-L^{\tau}(x,s)\|_{2p}& \leq & C_1\mid t-s\mid^{1-\tau}, \ \text{for} \ 0<\tau<1
\end{eqnarray}
and
\begin{eqnarray} \label{reg2}
&& \|L^{\tau}(y,t)-L^{\tau}(y,s)-L^{\tau}(x,t)+L^{\tau}(x,s)\|_{2p} \notag
\\
& \leq & C_2|y-x|^{\nu}\mid t-s\mid^{1-\tau(1+\nu)},
\end{eqnarray}
for any $0<\nu<\frac{1-\tau}{2\tau}$.\newline
\end{lemma}

Combining lemma \ref{LemReg} and a version of Kolmogorov's theorem in the
two-parameter H\"{o}lderian space $\mathcal{H}^{a,b}_{0}$, we obtain the
following result which will serve us in the sequel,

\begin{lemma}
Let $T> 0$. Then the following condition is satisfied almost surely: For all
$0 <\alpha<1-\tau$, there exists a constants $C_3> 0$ and $C_4>0$ such that
for all $(t, s)\in[0, T]^2$ and $x\in \mathbb{R}$ where $|x|\leq M$, we
have,
\begin{eqnarray} \label{eqRegT}
|L^\tau(x,t)-L^\tau(x,s)| & \leq & C_3|t-s|^{\alpha}.
\end{eqnarray}
For all $0 <\alpha<1-\tau$ and $0 <\beta<\frac{1-\tau}{2\tau}$, there exists
a constant $C_4> 0$ such that for all $(t, s)\in[0, T]^2$ and $(x, y)\in
\mathbb{R}^2$ where $|x|,|y|\leq M$, we have,
\begin{equation} \label{eq11}
|L^\tau(x,t)- L^\tau(y,t)-L^\tau(x,s) +L^\tau(y,s)| \leq
C_4|t-s|^{\alpha}|x-y|^{\beta}.
\end{equation}
\end{lemma}

\begin{remark}
By self-similarity of $X^{\tau}$, the local time $L^{\tau}$ verify the
so-called scaling property, in the sense,%
\begin{equation}
\{L^\tau(x,t), \ t\in \mathbb{R}^+\}\overset{(d)}{=}
\{\lambda^{\tau-1}L^\tau(x\lambda^\tau,\lambda t), \ t\in \mathbb{R}%
^+\},\ \forall \ \lambda>0. \label{eq1}
\end{equation}
\end{remark}

\begin{remark}
The previous constants $C_i$, $i\in\{1,2,3,4\}$ in the inequalities %
\eqref{reg1}, \eqref{reg2}, \eqref{eqRegT} and \eqref{eq11} and the
constants which will appear in the next results are random variables (unless
otherwise specified).
\end{remark}

Throughout the remainder of this paper, we will consider the jointly
continuous version of the local time of $X^{\tau}$ and we assume that this
process is the canonical process on the space $\mathcal{C}(\mathbb{R}^{+},%
\mathbb{R})$ of continuous functions.

We end this section with the following result which will serve us in the
next section:

\begin{lemma}
\label{lem1.1}The local times $L^\tau(x,t)$ of $X^{\tau}$ are additive
functionals with respect to time variable, that is, for all $t\geq s\geq0$
and for all $\omega \in \Omega$
\begin{equation}
L^\tau(x,t)(\omega)=L^\tau(x,s)(\omega)+L^\tau(x,t-s)(\theta_s(\omega))
\label{eq8}
\end{equation}
where for all $s\in \mathbb{R}^+$, $\theta_s$ is the shift operator defined
by $\theta_s(\omega)=\omega(t+s)$, i.e.
$X^\tau(t)\circ \theta_s=X^\tau(t+s).$
\end{lemma}

\begin{proof}
The local time $L^\tau$ has a representation of the form
\begin{eqnarray*}
L^\tau(x,t)=\int_{0}^{t}\delta(X^{\tau}(u)-x)du,
\end{eqnarray*}
where $\delta$ is the Dirac function at the point 0. If we put $%
X^\tau(u)=X^\tau_u$, then we have, for all $t\geq s\geq0$
\begin{eqnarray*}
L^\tau(x,t-s)(\theta_s(\omega))&=&\int_{0}^{t-s} \delta\left(
X^\tau_u(\theta_s(\omega))-x \right)du \\
&=&\int_{0}^{t-s} \delta\left(X^\tau_{u+s}(\omega)-x \right)du \\
&=&\int_{s}^{t} \delta\left(X^\tau_{u}(\omega)-x \right)du \\
\\
&=& \int_{0}^{t} \delta\left(X^\tau_{u}(\omega)-x \right)du -\int_{0}^{s}
\delta\left(X^\tau_{u}(\omega)-x \right)du \\
&=&L^\tau(x,t)(\omega)- L^\tau(x,s)(\omega .
\end{eqnarray*}
We can also find this property from equation \eqref{eq7}.
\hfill \end{proof}


\section{Strong approximation}

In this section, we will show the key ingredient which allows us to
establish the strong approximation of some additive functionals associated
to fBm and that associated to Riemann-Liouville process.

First, we need the following two propositions which are a consequence of the
scaling property \eqref{eq1} (for more applications of the scaling property
of fBm, we refer to Samorodnitsky and Taqqu \cite{ST}):

\begin{prop}
\label{prop2.1} For any $\lambda>0$, we have
\begin{equation}
\begin{aligned} & \ \ \ \left\{\sup_{(x,y)\in\mathbb{R}^2}
\left|L^\tau(x,t)-L^\tau(y,t)\right|, \ t\geq0, \
\mathbb{P}^0\right\}\overset{(d)}{=} \\
&\left\{\lambda^{-1+\tau}\sup_{(x,y)\in\mathbb{R}^2} \left|L^\tau(x,\lambda
t)-L^\tau(y,\lambda t)\right|, \ t\geq0, \ \mathbb{P}^0\right\}.\\ \label{eq5}
\end{aligned}
\end{equation}
\end{prop}
\newpage
\begin{prop}
\label{prop2.2}\textit{{For any $0<\gamma<\frac{1-\tau}{2\tau}$ and any $%
\lambda>0$, we have }}

\textit{%
\begin{equation}
\begin{aligned} & \ \ \ \left\{\sup_{0\leq s\leq t}\sup_{x\neq
y}\frac{\left|L^\tau(x,s)-L^\tau(y,s)\right|}{|x-y|^\gamma}, \ t\geq0, \
\mathbb{P}^0\right\}\overset{(d)}{=}\\ &\left\{\lambda^{-1+\tau(1+\gamma)}\sup_{0\leq
s\leq t}\sup_{x\neq y} \frac{\left|L^\tau(x,\lambda s)- L^\tau(y, \lambda
s)\right|}{|x-y|^\gamma}, \ t\geq0, \ \mathbb{P}^0\right\}.\\ \label{eq6}
\end{aligned}
\end{equation}
}
\end{prop}

Let us now prove the strong approximation version of the first-order limit
theorem for the additive functionals of the process $X^{\tau}$. To prove
this theorem, we still need of some results and a weak condition which
replace the spatial homogeneity property in the case of Markov process
defined as:

\begin{definition}
\label{lem2.211}\ A stochastic process $Y=\{Y(t), \ t\geq0\}$ is called
spatially homogeneous if its finite dimensional distributions do not change
with a shift in space.
\end{definition}

In the Brownian motion case $W=\{W(t), \ t\geq0\}$, the spatial homogeneous
property means that the transition probability (which has a density)
satisfies $p_t(x, y) = p_t(0, y-x)$ and then
\begin{equation*}
\mathbb{P}^{x+y}(W(t)\in A)=\mathbb{P}^x(W(t)\in A),
\end{equation*}
for all $x,y\in\mathbb{R}$ and any Borel subset $A$ of $\mathbb{R}$.
Unfortunately the process $X^{\tau}$ is not Markov and does not have such
property and then its distribution depends on the starting point. But the
associated tangent processes to its local time have the following property:

\begin{lemma}
\label{lem2.21} For all $t\in(0,+\infty)$ and all $\nu>0$, the stochastic
process $Z=\{Z(t), \ t\geq0\}$ defined by
\begin{equation*}
Z(t):=\sup_{x\neq y}\frac{\left|L^\tau(x,t)-L^\tau(y,t)\right|}{|x-y|^{\nu}},
\end{equation*}
satisfies the following equality in law
\begin{equation*}
\mathbb{P}^z(Z(t)\in A)=\mathbb{P}^0(Z(t)\in A), \ z\in\mathbb{R}
\end{equation*}
for all Borel subset $A$ of $\mathbb{R}$.
\end{lemma}


\begin{proof}
Let $z\in \mathbb{R}$ and consider $\sigma _{z}$ the translation mapping of $%
\Omega \longrightarrow \Omega $ defined by $\sigma _{z}(\omega )=\omega
(\cdot )+z$. This means that the translation operator $\sigma _{z}$ maps a
path $t\mapsto \omega (t)$ to the new path $t\mapsto \omega (t)+z$. So,
indeed if we put $X_{t}^{\tau }=X^{\tau }(t)$ we have, $X_{t}^{\tau }(\sigma
_{z}(\omega ))=X_{t}^{\tau }(\omega )+z$ and the translated process $%
\{X^{\tau }(t)+z),\ t\geq 0\}$, under the law $\mathbb{P}^{x}$, has the same
distribution as $(X^{\tau }(t),\ t\geq 0)$ under $\mathbb{P}^{x+y}$. Then for all $%
x\in \mathbb{R}$ and any Borel subset $A$ of $\mathbb{R}$
\begin{equation*}
\mathbb{P}^{x+y}(X^{\tau }(t)\in A)=\mathbb{P}^{x}(X^{\tau }(t)+y\in A).
\end{equation*}%
Hence, from (\ref{eq08}), we have for any Borel subset $A$ of $\mathbb{R}$
{\small \begin{eqnarray*}
&&\mathbb{P}^{z}\left( Z(t)\in A\right) \\
&=&\mathbb{P}^{z}\left( \sup_{x\neq y}\frac{\left\vert \lim_{\varepsilon \rightarrow
0^{+}}\frac{1}{2\varepsilon }\left[ \int_{0}^{t}\mathbb{I}_{\mid X^{\tau
}(s)-x\mid \leq \varepsilon }ds-\int_{0}^{t}\mathbb{I}_{\mid X^{\tau
}(s)-y\mid \leq \varepsilon }ds\right] \right\vert }{|x-y|^{\nu }}\in
A\right) \\
&=&\mathbb{P}^{0}\left( \sup_{x\neq y}\frac{\left\vert \lim_{\varepsilon \rightarrow
0^{+}}\frac{1}{2\varepsilon }\left[ \int_{0}^{t}\mathbb{I}_{\mid X^{\tau
}(s)-x-z\mid \leq \varepsilon }ds-\int_{0}^{t}\mathbb{I}_{\mid X^{\tau
}(s)-y-z\mid \leq \varepsilon }ds\right] \right\vert }{|x-y|^{\nu }}\in
A\right) \\
&=&\mathbb{P}^{0}\left( \sup_{x\neq y}\frac{\left\vert L^{\tau }(x+z,t)-L^{\tau
}(y+z,t)\right\vert }{|x-y|^{\nu }}\in A\right) \\
&=&\mathbb{P}^{0}\left( \sup_{x\neq y}\frac{\left\vert L^{\tau }(x,t)-L^{\tau
}(y,t)\right\vert }{|x-y|^{\nu }}\in A\right) \\
&=&\mathbb{P}^{0}(Z(t)\in A).
\end{eqnarray*}}
This ends the proof \hfill \end{proof}


\begin{lemma}
\label{lem2.2} Let $y\in \mathbb{R}$. Then for all $t\in (0,+\infty )$ and
all $0<\nu <\frac{1-\tau }{2\tau }$, the process $Y=\{Y(t),\ t\geq 0\}$
defined by
\begin{equation}
Y(t):=\sup_{0\leq s\leq t}\sup_{x\neq y}\frac{\left\vert L^{\tau
}(x,s)-L^{\tau }(y,s)\right\vert }{|x-y|^{\nu }} \label{Y}
\end{equation}%
is increasing and continuous.
\end{lemma}

\begin{proof}
As $Y(0)=0$, we can see easily that $\{Y(t),\ t>0\}$ is an increasing
process. On the other hand, for all $0\leq u\leq s\leq t$, we have
\begin{eqnarray*}
&&\sup_{x\neq y}\frac{\left\vert L^{\tau }(x,s)-L^{\tau }(y,s)\right\vert }{%
|x-y|^{\nu }} \\
&=&\sup_{x\neq y}\frac{\left\vert L^{\tau }(x,u)+L^{\tau }(x,s-u)\circ
\theta _{u}-L^{\tau }(y,u)-L^{\tau }(y,s-u)\circ \theta _{u}\right\vert }{%
|x-y|^{\nu }} \\
&\leq &\sup_{x\neq y}\frac{\left\vert L^{\tau }(x,u)-L^{\tau
}(y,u)\right\vert }{|x-y|^{\nu }} \\ && +\sup_{x\neq y}\frac{\left\vert L^{\tau
}(x,s-u)\circ \theta _{u}-L^{\tau }(y,s-u)\circ \theta _{u}\right\vert }{%
|x-y|^{\nu }} \\
&\leq &\sup_{0\leq v\leq u}\sup_{x\neq y}\frac{\left\vert L^{\tau
}(x,u)-L^{\tau }(y,u)\right\vert }{|x-y|^{\nu }} \\ && +\sup_{x\neq y}\frac{%
\left\vert L^{\tau }(x,s-u)\circ \theta _{u}-L^{\tau }(y,s-u)\circ \theta
_{u}\right\vert }{|x-y|^{\nu }} \\
&\leq &Y(u)+\sup_{x\neq y}\frac{\left\vert L^{\tau }(x,s-u)-L^{\tau
}(y,s-u)\right\vert }{|x-y|^{\nu }}\circ \theta _{u}.
\end{eqnarray*}%
Therefore
\begin{eqnarray*}
&&\sup_{0\leq s\leq t}\sup_{x\neq y}\frac{\left\vert L^{\tau }(x,s)-L^{\tau
}(y,s)\right\vert }{|x-y|^{\nu }} \\
&\leq &Y(u)+\sup_{0\leq s\leq t}\sup_{x\neq y}\frac{\left\vert L^{\tau
}(x,s-u)-L^{\tau }(y,s-u)\right\vert }{|x-y|^{\nu }}\circ \theta _{u} \\
&\leq &Y(u)+\sup_{0\leq s-u\leq t-u}\sup_{x\neq y}\frac{\left\vert L^{\tau
}(x,s-u)-L^{\tau }(y,s-u)\right\vert }{|x-y|^{\nu }}\circ \theta _{u} \\
&=&Y(u)+\sup_{0\leq v\leq t-u}\sup_{x\neq y}\frac{\left\vert L^{\tau
}(x,v)-L^{\tau }(y,v)\right\vert }{|x-y|^{\nu }}\circ \theta _{u}.
\end{eqnarray*}%
Thus we have
\begin{equation*}
Y(t)\leq Y(u)+Y(t-u)\circ \theta _{u}.
\end{equation*}%
By the inequality (\ref{eq11}), we deduce that for any $0<\beta <\frac{%
1-\tau }{2\tau }$,
\begin{equation*}
0\leq Y(t)-Y(u)\leq Y(t-u)\circ \theta _{u}\leq C_{4}\sup_{u \leq s\leq
t}|s-u|^{\beta }\leq C_{4}|t-u|^{\beta }.
\end{equation*}%
This means that $\{Y(t),\ t>0\}$ is H\"{o}lder continuous.
\hfill \end{proof}

\begin{remark}
We may reason otherwise to prove the continuity of $Y$ by using \cite[Lemma
7.30, p. 207]{Mo} : the upper bound of every continuous function $%
f:[0,t]\rightarrow \mathbb{R}$ $(t>0)$, can be written as follows
\begin{equation*}
\sup_{0\leq s\leq t}f(s)=\lim_{a\rightarrow+\infty}\frac{1}{a}%
\ln\left(\int_{0}^{t}\exp(af(v))dv\right).
\end{equation*}
Using this lemma with
\begin{equation*}
f(s)=\sup_{x\neq y}\frac{\left|L^\tau(x,s)-L^\tau(y,s)\right|}{|x-y|^{\gamma}%
},
\end{equation*}
which is continuous, we get the continuity of $t\mapsto Y(t)$.
\end{remark}

\begin{remark}
In the precedent lemma, we have shown the H\"olderian continuity of the
process $Y$, but the continuity of $Y$ is enough for us to prove the
following corollary:
\end{remark}

\begin{cor}
\label{cor2.2} For all $t\geq 0$, $0<\nu <\frac{1-\tau }{2\tau }$ and $p\geq 1$, there exists a non random constant $0<C_{5}<+\infty $, such
that
\begin{equation*}
\left\Vert Y(t)\right\Vert _{p}\leq C_{5}t^{1-\tau (1+\nu )}
\end{equation*}
\end{cor}

\begin{proof}
Let $\{Y(t),\ t\geq 0\}$ be as defined in (\ref{Y}), hence by the Lemma \ref{lem2.21} and the equality (in law) (\ref{eq6}), we have
{\small \begin{eqnarray*}
\lim_{z\rightarrow+\infty}\sup_{\substack{ \lambda>0, \\ u\in\mathbb{R}}}%
\mathbb{P}^u\left(Y(\lambda)>\lambda^{1-\tau(1+\nu)}z \right)
&=&\lim_{z\rightarrow+\infty}\sup_{\lambda>0}\mathbb{P}^0\left(Y(\lambda)>\lambda^{1-%
\tau(1+\nu)}z \right) \\
&=&\lim_{z\rightarrow+\infty}\mathbb{P}^0\left(Y(1)>z \right) \\
&=&0,
\end{eqnarray*}}
where we have used in the last estimation the inequality (\ref{eq11}) and
the fact that $L^{\tau }(x,0)=0$, implies the random variable $Y(1)$ is
finite a.s. Therefore by the Lemma \ref{lem2.2} and according to the \cite[%
Lemma 4.5, p. 326]{FG}, applied to $\{Y(t),\ t\geq 0\}$ and $q=\frac{1}{%
1-\tau (1+\nu )}$ we deduce the result of our Corollary \ref{cor2.2}.
\hfill \end{proof}

\begin{cor}
\label{cor2.1}For any $0<\tau<1$, $0 < \nu <\frac{1-\tau}{2\tau} $ and $\varepsilon>0$,
for $t$ sufficiently large, we have
\begin{eqnarray}
\sup_{x\neq y}\frac{\left|L^\tau(x,t)-L^\tau(y,t)\right|}{|x-y|^{\nu}}%
&=&o\left( t^{1-\tau(1+\nu)+\varepsilon}\right) \ \text{\ a.s. }
\label{eq14}
\end{eqnarray}
\end{cor}

\begin{proof}
Using Tchebychev's inequality and the precedent corollary with $p = \frac{2}{%
\varepsilon}$, for any $n \geq 1$, there is $C(\nu,\varepsilon)$ such that

\begin{equation*}
\mathbb{P}^{0}\left( Y(n)>n^{1-\tau (1+\nu )+\varepsilon }\right) \leq \frac{C(\nu
,\varepsilon )}{n^{2}},
\end{equation*}
and by the Borel-Cantelli lemma, we deduce as $n\longrightarrow +\infty $,
\begin{equation*}
Y(n)=\mathcal{O}\left( t^{1-\tau (1+\nu )+\varepsilon }\right) \ \text{\
a.s. }
\end{equation*}%
Since the process $\{Y(t),\ t\geq 0\}$ is increasing and $\varepsilon $ can
be arbitrarily small, we have proved the desired corollary. \hfill \end{proof}
\newpage
We will also need the following estimation:
\begin{lemma}
\label{lem2.3} For any $\varepsilon > 0$, when $t$ goes to infinity, we have
\begin{eqnarray}
\sup_{(x,y)\in\mathbb{R}^{2}}\left|L^\tau(x,t)-L^\tau(y,t)\right|&=& o\left(
t^{1-\tau+\varepsilon}\right) \ \text{\ a.s. } \label{eq15}
\end{eqnarray}
\end{lemma}

\begin{proof}
The process $\{K(t),\ t\geq 0\}$ defined by
\begin{equation*}
K(t):=\sup_{0\leq s\leq t}\sup_{(x,y)\in\mathbb{R}^2}\left|L^\tau(x,s)-L^%
\tau(y,s)\right|,
\end{equation*}
is increasing and continuous. On the other hand, we can proceed as Lemma \ref{lem2.21} to prove that the process $\{K(t), \ t\geq0\}$ 
satisfies the following equality in law
\begin{equation*}
\mathbb{P}^z(K(t)\in A)=\mathbb{P}^0(K(t)\in A), \ z\in\mathbb{R}
\end{equation*}
for all Borel subset $A$ of $\mathbb{R}$.
Furthermore, by the equality in law (\ref{eq5}), we have
\begin{equation*}
K(t)\overset{(d)}{=}\lambda^{-1+\tau}K(\lambda t),
\end{equation*}
which gives in particular $K(\lambda)\overset{(d)}{=}\lambda^{1-\tau}K(1)$.
Now we apply \cite[Lemma 4.5, p. 326]{FG} to the process $\{K(t),\ t\geq 0\}
$ and $q=\frac{1}{1-\tau}$, by the same arguments used in the proof of the
precedent corollary we obtain
\begin{equation*}
\sup_{(x,y)\in\mathbb{R}^{2}}\left|L^\tau(x,t)-L^\tau(y,t)\right|= o\left(
t^{1-\tau+\varepsilon}\right) \ \text{\ a.s. }
\end{equation*}
This the desired result. \hfill \end{proof}

Now we are ready to prove the following result of strong approximation:

\begin{theorem}
\label{thm2.2} Let $f$ be a Borel function on $\mathbb{R}$ such that $%
\overline{f}\neq0$ and
\begin{eqnarray}
\int_{\mathbb{R}}|x|^k |f(x)|dx &<& +\infty, \label{eq16a}
\end{eqnarray}
for some $k>0$. Then, for all sufficiently small $\varepsilon>0$, when $t$
goes to infinity, we have
\begin{eqnarray}
\int_{0}^{t}f\left(X^{\tau}(s)\right)ds&=& \overline{f}L^\tau(0,t)+o\left(
t^{1-\tau-\varepsilon} \right)\ \ \text{\ a.s. } \label{eq16}
\end{eqnarray}
\end{theorem}

\begin{proof}
Without loss of generality, we suppose that $f$ is a real positive function, for the general case we use the decomposition $f=f^{+}-f^{-}$. \\
Let
\begin{equation*}
J(t):=\int_{0}^{t}f\left( X^{\tau }(s)\right) ds-\overline{f}L^{\tau }(0,t).
\end{equation*}%
By the occupation density formula (\ref{eq3}), we have
\begin{equation*}
I(t)=\int_{\mathbb{R}}f\left( x\right) \left( L^{\tau }(x,t)-L^{\tau
}(0,t)\right) dx=J_{1}(t)+J_{2}(t),
\end{equation*}%
where
\begin{equation*}
J_{1}(t):=\int_{\{|x|>t^{a}\}}f\left( x\right) \left( L^{\tau }(x,t)-L^{\tau
}(0,t)\right) dx,
\end{equation*}%
and
\begin{equation*}
J_{2}(t):=\int_{\{|x|\leq t^{a}\}}f\left( x\right) \left( L^{\tau
}(x,t)-L^{\tau }(0,t)\right) dx,
\end{equation*}%
for some $0<a<1$. Let us deal with the first term $J_{1}(t)$. By the
hypothesis (\ref{eq16a}) and the estimation (\ref{eq15}), we have for some
sufficiently small $\varepsilon >0$,
\begin{eqnarray*}
\left\vert J_{1}(t)\right\vert &\leq& \sup_{\{|x|>t^{a}\}}\left\vert L^{\tau
}(x,t)-L^{\tau }(0,t)\right\vert
\int_{\{|x|>t^{a}\}}|x|^{-k}|x|^{k}\left\vert f(x)\right\vert dx \\
&=&o\left( t^{1-\tau -ak+\varepsilon }\right) \text{ \ a.s. }
\end{eqnarray*}%
On the other hand, by using the estimation (\ref{eq14}) and the fact that $%
f $ is integrable, we deduce that, for any $0<\nu <\frac{1-\tau }{2\tau }$
and $\varepsilon >0$,
\begin{eqnarray*}
\left\vert J_{2}(t)\right\vert &\leq& t^{a\nu }\sup_{\{|x|\leq t^{a}\}}\frac{%
\left\vert L^{\tau }(x,t)-L^{\tau }(0,t)\right\vert }{|x|^{\nu }}%
\int_{\{|x|\leq t^{a}\}}\left\vert f(x)\right\vert dx \\
&=&o\left( t^{1-\tau (1+\nu )+a\nu +\varepsilon }\right) \text{ \ a.s. }
\end{eqnarray*}%
It follows that
\begin{equation*}
\left\vert J(t)\right\vert =o\left( t^{1-\tau -ak+\varepsilon }\right)
+o\left( t^{1-\tau (1+\nu )+a\nu +\varepsilon }\right) \text{ \ a.s. }
\end{equation*}%
Choosing $a=\frac{\nu \tau }{\nu +k}$. It is clear that $0<a<1$ and then
\begin{equation*}
\left\vert J(t)\right\vert =o\left( t^{1-\tau -ka+\varepsilon }\right) \text{
\ a.s. }
\end{equation*}%
Then for all sufficiently small $\varepsilon >0$ and $t$ sufficiently large,
we deduce that
\begin{equation*}
J(t)=o\left( t^{1-\tau -\varepsilon }\right) \text{ \ a.s. }
\end{equation*}%
This completes the proof of the theorem.
\hfill \end{proof}


\section{Law of Iterated Logarithm for some additive functionals of the
process $X^{\protect\tau}$}

In this section and as application of strong approximation of the Theorem %
\ref{thm2.2}, we obtain the Law of Iterated Logarithm (LIL) of the process
\begin{equation*}
\left\{ \int_{0}^{t}f\left( X^{\tau }(s)\right) ds,\ t\geq 0\right\} .
\end{equation*}%
Let $\ell (x,t)$ be the local time of the standard Brownian motion \mbox{$(W(t),\ t\geq 0)$.} Kesten \cite{K} showed that $\ell (x,t)$ satisfies the following LIL,
\begin{equation*}
\limsup_{t\rightarrow +\infty }\frac{\ell (0,t)}{\sqrt{2t\log \log t}}%
=\limsup_{t\rightarrow +\infty }\sup_{x\in \mathbb{R}}\frac{\ell (x,t)}{%
\sqrt{2t\log \log t}}=1 \ \text{\ a.s. }
\end{equation*}%
For the LIL of the local time $L^{\alpha }(x,t)$ of the so-called symmetric $%
\alpha $-stable process $X^{\alpha }$ was given by Donsker and Varadhan \cite%
{DV}, as follows
\begin{equation*}
\limsup_{t\rightarrow +\infty }\frac{L^{\alpha }(0,t)}{t^{\frac{\alpha -1}{%
\alpha }}\left( \log \log t\right) ^{\frac{1}{\alpha }}}=\limsup_{t%
\rightarrow +\infty }\sup_{x\in \mathbb{R}}\frac{L^{\alpha }(x,t)}{t^{\frac{%
\alpha -1}{\alpha }}\left( \log \log t\right) ^{\frac{1}{\alpha }}}=d(\alpha
),\ \ \text{\ a.s. }
\end{equation*}%
where
\begin{equation*}
d(\alpha ):=\frac{\Gamma \left( \frac{1}{\alpha }\right) \Gamma \left( 1-%
\frac{1}{\alpha }\right) }{\pi (\alpha -1)^{1-\frac{1}{\alpha }}}.
\end{equation*}%
Ait Ouahra and Sghir \cite{AS} extended this property for a class of some
additive functionals of $X^{\alpha }$, more precisely, they studied the LIL
of the process
\begin{equation*}
\left\{ \int_{0}^{t}f\left( X_{s}^{\alpha }\right) ds,\ t\geq 0\right\} ;
\end{equation*}%
that is,
\begin{equation*}
\limsup_{t\rightarrow +\infty }\frac{\int_{0}^{t}f\left( X_{s}^{\alpha
}\right) ds}{t^{\frac{\alpha -1}{\alpha }}\left( \log \log t\right) ^{\frac{1%
}{\alpha }}}=\overline{f}d\left( \alpha \right) ,\ \ \text{\ a.s. }
\end{equation*}%
where $f$ is a Borel function on $\mathbb{R}$ such that $\overline{f}\neq 0$
and $\int_{\mathbb{R}}|x|^{k}|f(x)|dx$ is finite for some $k>0$.

Chen \emph{et al.} \cite{CLRS} treated the LIL of the local time $L^{\tau
}(x,t)$ (in both cases: fBm and Riemann-Liouville process) and given by the
following almost sure equality
\begin{equation}
\limsup_{t\rightarrow +\infty }\frac{L^{\tau }(0,t)}{t^{1-\tau }\left( \log
\log t\right) ^{\tau }}=c_{\tau }\left( \theta (\tau )\right) ^{-\tau },\ \
\text{\ a.s. } \label{eq4.32}
\end{equation}%
with
\begin{equation*}
\left( \frac{\pi \delta _{\tau }^{2}}{\tau }\right) ^{1/2\tau }\theta
_{0}(\tau )\leq \theta (\tau )\leq \left( 2\pi \right) ^{1/2\tau }\theta
_{0}(\tau )
\end{equation*}%
where
\begin{equation}
\delta _{\tau }=\sqrt{2\tau }2^{\tau }\beta (1-\tau ,\tau +1/2)^{-1/2}, \
\theta _{0}(\tau )=\tau \left( \frac{(1-\tau )^{1-\tau }}{\Gamma (1-\tau )}%
\right) ^{1/\tau } \label{eq4.33}
\end{equation}%
and $c_{\tau }=\delta _{\tau }$ for Riemann-Liouville process case, and $%
c_{\tau }=1$ for fBm case. \newline
Notice that the inequalities in (\ref{eq4.32}) become equalities in the
Brownian motion case $\tau =\frac{1}{2}$.\newline
The Theorem \ref{thm2.2} allow us to deduce the following LIL:

\begin{theorem}
\label{thm3.1}Let $f$ be a Borel function on $\mathbb{R}$ such that $%
\overline{f}\neq0$ and
\begin{equation*}
\int_{\mathbb{R}}|x|^k |f(x)|dx<+\infty,
\end{equation*}
for some $k>0$. Then it holds that
\begin{eqnarray}
\limsup_{t\rightarrow +\infty}\frac{\int_0^t f \left(X^{\tau}(s) \right)ds}{%
t^{1-\tau} \left(\log\log t\right)^\tau}&=&\overline{f}c_{\tau}\left(\theta(%
\tau)\right)^{-\tau} \ \text{\ a.s. } \label{eq4.4}
\end{eqnarray}
where $c_{\tau}$ is the constant appeared in (\ref{eq4.32}).
\end{theorem}

\begin{proof}
According to the theorem \ref{thm2.2}, for all sufficiently small $%
\varepsilon>0$, when $t$ goes to infinity, we have
\begin{equation*}
\int_{0}^{t}f\left(X^{\tau}(s)\right)ds= \overline{f}L^\tau(0,t)+o\left(
t^{1-\tau-\varepsilon} \right), \ \ \text{\ a.s. }
\end{equation*}
and the fact that
\begin{equation*}
\limsup_{t\rightarrow\infty}\frac{L^\tau(0,t)}{t^{1-\tau}\left(\log\log t
\right)^{\tau}}=c_{\tau}\left(\theta(\tau)\right)^{-\tau}, \ \ \ a.s.
\end{equation*}
we obtain almost surely
\begin{eqnarray*}
\limsup_{t\rightarrow\infty}\frac{\int_0^t f \left(X^{\tau}(s) \right)ds}{%
t^{1-\tau} \left(\log\log t\right)^\tau}&=&\overline{f} \limsup_{t%
\rightarrow\infty}\frac{L^\tau(0,t)}{t^{1-\tau}\left(\log\log t
\right)^{\tau}} \\
&=&c_{\tau}\overline{f}.\left(\theta(\tau)\right)^{-\tau}.
\end{eqnarray*}
This achieves the goal. \hfill \end{proof}

\begin{remark}
These results can be extended to other self-similar Gaussian processes
namely: the bi-fractional Brownian motion and sub-fractional Brownian
motion, by giving the limit theorem and the strong approximation versions
and that of the associated LIL of the same additive functionals studied in
this paper.
\end{remark}


\noindent Mohamed Ait Ouahra \\ Department of Mathematics, Faculty of Sciences, Mohamed First University B.P. 717 - 60000, Oujda, Morocco, \\ \emph{E-mail address}: m.aitouahra@ump.ac.ma \\ \ \\
Abderrahim Aslimani \\ Department of Mathematics Pluridisciplinary
Faculty, Mohamed First University 300-62700,
Selouane, Nador Morocco, \\ \emph{E-mail address}: a.slimani@ump.ac.ma \\ \ \\
Mhamed Eddahbi \\ Department of Mathematics, College of Sciences,
King Saud University, PO. Box 2455, Z.C. 11451 Riyadh, Saudi Arabia \\ \emph{E-mail address}: meddahbi@ksu.edu.sa \\ \ \\
Mohamed Mellouk \\ Univesit\'e Paris Cit\'e, CNRS, MAP5, F-75006, Paris, France, \\ \emph{E-mail address}: mohamed.mellouk@parisdescartes.fr
\end{document}